\documentclass{amsart}

\usepackage[T1]{fontenc}
\usepackage[utf8]{inputenc}
\usepackage{lmodern}
\usepackage{appendix}
\usepackage{amsmath,amssymb,amsthm,mathtools}
\usepackage{hyperref}
\usepackage[nameinlink,noabbrev]{cleveref}
\usepackage{graphicx}

\newtheorem{theorem}{Theorem}[]
\newtheorem*{theorem*}{Theorem}
\newtheorem{corollary}{Corollary}
\newtheorem{lemma}{Lemma}
\newtheorem{conj}{Conjecture}

\newtheorem{claim}{Claim}

\usepackage{float}

\newcommand{\G}{G_{27}}
\newcommand{\Gp}{G_{29}}

\title[A unit-distance graph in the plane with independence ratio below 1/4]
{A unit-distance graph in the plane with independence ratio below 1/4}
\author[Á. Dúcz]{Ákos Dúcz}
\author[D. Varga]{Dániel Varga}
\date{June 2026}

\begin{document}

\begin{abstract}
We prove that there exists a finite unit-distance graph in the plane with independence ratio strictly smaller than $1/4$, answering a question of Erdős. Our proof closely follows the framework of Matolcsi, Ruzsa, Varga, and Zsámboki~\cite{MRVZ}, based on the geometric fractional chromatic number, but adds a carefully chosen two-vertex augmentation that pushes their $27$-vertex construction from geometric fractional chromatic number $4$ to a value strictly larger than $4$. This disproves their Conjecture 1, and implies that the fractional chromatic number of the plane is strictly larger than $4$. The proof can be made fully constructive, but the resulting finite graph has an enormous number of vertices.

\end{abstract}

\vspace*{-1.5cm}

\maketitle

\section{Introduction}

In a 1987 problem list, Erdős posed the following question. Let $f(n)$ be the largest integer such that every set of $n$ points in the plane contains $f(n)$ points no two of which are at distance $1$. Equivalently, $f(n)$ is the minimum possible independence number of a unit-distance graph arising from $n$ points in the plane. Erdős asked to determine $f(n)$ as accurately as possible and, in particular, whether $f(n)\ge n/4$ for all $n$ \cite{Er87}. In this note, we answer this question in the negative:

\begin{theorem}\label{thm1}
There exists a unit-distance graph $G$ in the plane such that
\[
    \frac{\alpha(G)}{|V(G)|}<\frac14.
\]
\end{theorem}

We only prove the existence of this graph here. In principle, the construction can be made fully explicit, but the resulting graph is astronomically large.

Our result falsifies Conjecture 1 of Matolcsi et al. \cite{MRVZ}, which states that for any finite unit-distance graph $G$, its fractional chromatic number $\chi_f(G)$ is strictly smaller than 4:

\begin{corollary}
\[
    \chi_f(\mathbb R^2) > 4.
\]
\end{corollary}

As immediate consequences, we recover de Grey's lower bound \cite{deGrey},
\begin{corollary}
\[
    \chi(\mathbb R^2) \geq 5,
\]
\end{corollary}

and obtain a simpler proof of another conjecture of Erdős, concerning the density of measurable $1$-avoiding sets in the plane, first proved by Ambrus et al. \cite{ACMVZ} using Fourier-analytic tools:

\begin{corollary}
\[
    m_1(\mathbb{R}^2) < 1/4.
\]
\end{corollary}

We closely follow the framework of Matolcsi, Ruzsa, Varga, and Zsámboki \cite{MRVZ}, which is based on the geometric fractional chromatic number $\chi_{gf}(G)$. Their method consists of a three-step strategy for constructing unit-distance graphs with small independence ratio. In the present paper, we modify only the first step, and then use the remaining two steps unchanged.

\begin{enumerate}
    \item They construct a specific 27-vertex point configuration $G_{27}$ with $\chi_{gf}(G_{27}) = 4$.
    \item Exploiting the amenability of the group of Euclidean transformations in dimension 2, they apply a blow-up procedure placing many congruent copies of $\G$ to construct unit-distance graphs $G'$ with $\chi_f(G')$ arbitrarily close to $\chi_{gf}(\G)$.
    \item Starting from some $G'$, they apply a second blow-up procedure, based on large discrete cubes in the additive lattice generated by $V(G')$, to construct unit-distance graphs with independence ratio arbitrarily close to $1/\chi_f(G')$.
\end{enumerate}

In this note, we revisit their $27$-vertex point configuration $G_{27}$ and show that adding two suitably chosen points is enough to push its geometric fractional chromatic number strictly above $4$. As we will see, such augmentations are exceedingly rare. Most of the note is devoted to finding one, and proving that it has the desired effect. To this end, we characterize all optimal geometric fractional colorings of $G_{27}$: they form the convex hull of only $23$ extremal colorings and have affine dimension $11$. The small number of extremal colorings is surprising, since the underlying linear program has $182304$ variables and $16855$ constraints. Moreover, these extremal colorings have small denominators; each can be realized as a $k$-fold coloring with $k\le 22$. 

Our search for augmentations of $G_{27}$ is guided by this characterization of the optimal face. Since the extremal optimal geometric fractional colorings are few in number and have small denominators, we can test candidate vertices by solving finite constraint satisfaction problems. For each candidate vertex, we ask which colorings can be extended to the new vertex. A vertex that obstructs many colorings is then regarded as a promising candidate. We then identify two of the most promising candidates, together obstructing all our colorings. Adding both of them to $G_{27}$ brings its geometric fractional chromatic number above 4, which can be certified by a dual solution to the appropriate linear program. We describe this procedure in detail in the main text. 

\section{Previous work}

\subsection{The fractional chromatic number of the plane}

Originally introduced as a linear programming relaxation of the famous Hadwiger--Nelson problem \cite{Gardner60,Soifer}, a major motivation for studying $\chi_f(\mathbb R^2)$ came from a measurable analogue of the problem stated above, also asked by Erdős. Let $m_1(\mathbb R^2)$ denote the supremum of the upper densities of measurable sets $A\subseteq\mathbb R^2$ containing no two points at distance~$1$. Erdős conjectured that $m_1(\mathbb R^2) < 1/4$ \cite{Er85}.

It is easy to see the following:

\begin{claim}\label{m1_vs_chi_f}
For any unit-distance graph $G$, $m_1(\mathbb R^2) \le 1 / \chi_f(G)$. 
\end{claim}

This claim is the weighted version of the standard averaging argument giving $m_1(\mathbb R^2)\le \alpha(G)/|V(G)|$. Indeed, the reciprocal $1/\chi_f(G)$ is the minimum, over all nonzero nonnegative vertex weightings $w$ of $G$, of the ratio $\max\{w(I): I\in\mathcal I(G)\}/w(V(G))$. The Croft construction \cite{Croft} gives $m_1(\mathbb R^2) \ge \delta_{\text{Croft}} = 0.22936\dots$ Together with Claim~\ref{m1_vs_chi_f}, this implies $\chi_f(\mathbb R^2) \leq 1/\delta_{\text{Croft}} < 4.36$.

A natural way of proving $m_1(\mathbb R^2)<1/4$ is therefore to find finite unit-distance graphs $G$ with $\chi_f(G)>4$. This approach produced a long sequence of lower bounds for $\chi_f(\mathbb R^2)$: The Moser spindle already gives $\chi_f(\mathbb R^2)\ge 7/2$, and this bound was subsequently improved by Mahan to $144/41$, by Fisher and Ullman to $32/9$, and by Cranston and Rabern to $76/21$ \cite{Mahan,FisherUllman,CranstonRabern}. Further progress was obtained by Bellitto, Pêcher, and Sédillot, who found a $607$-vertex graph giving the lower bound $3.899$ \cite{BellittoPecherSedillot}. Unpublished computations of Parts pushed this finite-graph approach still closer to $4$, reaching approximately $3.9898$ on a graph with $1057$ vertices \cite{Parts}. Despite these increasingly large constructions, no finite unit-distance graph with fractional chromatic number greater than $4$ was known before the current work.

The conjecture $m_1(\mathbb R^2)<1/4$ was eventually settled by Ambrus et al.~\cite{ACMVZ}, not by establishing $\chi_f(\mathbb R^2)>4$, but by combining the geometric fractional chromatic number method with Fourier analytic techniques.

Later, Matolcsi et al. \cite{MRVZ} gave a Fourier-free, geometric fractional chromatic number-based bound reaching the critical threshold: they proved $\chi_f(\mathbb R^2)\ge 4$. This gives only $m_1(\mathbb R^2)\le 1/4$, whereas Erdős asked for the strict inequality. The present note crosses this threshold by producing a finite unit-distance graph with independence ratio $\alpha(G)/|V(G)|$ strictly smaller than $1/4$, therefore providing a Fourier-free proof of the conjecture.

The task of estimating $f(n)/n$ accurately is still wide open. In fact, we conjecture

\begin{conj}
$f(n)/n = m_1(\mathbb R^2) + o(1)$
\end{conj}
and in particular that $m_1(\mathbb R^2) = \delta_{\text{Croft}}$.

\subsection{The Moser lattice}

Since our augmentation of $G_{27}$ falsifies \cite[Conjecture 1]{MRVZ} claiming that the independence ratio of a unit distance graph in the plane is always larger than $1/4$, it is worth revisiting the line of reasoning that led to the conjecture. The Moser lattice is the additive group generated by the elements of the Moser spindle in its standard embedding into the complex plane. The relevance of the Moser lattice first became apparent in the work of Ambrus, Csiszárik, Matolcsi, Varga, and Zsámboki \cite{ACMVZ}: although their search for unit-distance graphs was carried out in a much larger ambient ring in the complex plane, the graphs giving the best bounds turned out to be subsets of the Moser lattice.

Motivated by this observation, Matolcsi, Ruzsa, Varga, and Zsámboki \cite{MRVZ} pursued a search restricted to the Moser lattice, which produced the graph $G_{27}$ with $\chi_{gf}(G_{27})=4$. The lattice-restricted strategy then stopped improving at this value, a phenomenon later explained by the first author of the present note, who showed that the entire Moser lattice is geometrically integer $4$-colorable \cite{DuczMoserLattice}. In unpublished experiments while preparing \cite{MRVZ}, Matolcsi et al. also carried out searches not restricted to the lattice, in two variants: one starting from the Moser spindle, which failed to reach $\chi_{gf}=4$, and another starting from the already constructed graph $G_{27}$, which failed to improve on the value 4.

The above pieces of evidence made \cite[Conjecture 1]{MRVZ} a natural conjecture. Our current results help explain why the counterexample was difficult to find: the points that improve the geometric fractional chromatic number of $G_{27}$ appear to be extremely rare. As we discuss below, the successful augmentation was found only after using the structure of the optimal face of the geometric fractional coloring polytope to guide the search.

\section{Notation and preliminaries}

Throughout the note, a unit-distance graph means a finite point configuration in the plane together with the graph whose edges join pairs of points at distance $1$. Thus the embedding into the plane is part of the information, and when we talk about the distance of two vertices, we always mean Euclidean distance rather than graph theoretic distance. If $G$ is such a graph, we write $V(G)$ for its vertex set and $\mathcal I(G)$ for the set of all independent sets of $G$.

For a finite graph $G$, let $\alpha(G)$ denote its independence number. The ratio
$\alpha(G)/|V(G)|$ is called the independence ratio of $G$.

A fractional coloring of $G$ is a function
\[
    \gamma:\mathcal I(G)\to \mathbb R_{\ge 0}
\]
such that
\[
    \sum_{\substack{I\in\mathcal I(G)\\ v\in I}} \gamma(I) \ge 1
    \qquad \text{for every } v\in V(G).
\]
The weight of $\gamma$ is
\[
    |\gamma|=\sum_{I\in\mathcal I(G)} \gamma(I),
\]
and the fractional chromatic number $\chi_f(G)$ is the minimum possible value
of $|\gamma|$ over all fractional colorings of $G$.

We now recall the geometric variant introduced by \cite{ACMVZ}. If
$S\subseteq V(G)$, write
\[
    \overline\gamma(S)
    =
    \sum_{\substack{I\in\mathcal I(G)\\ S\subseteq I}} \gamma(I).
\]
Thus $\overline\gamma(S)$ is the total weight of colors assigned simultaneously to all
vertices of $S$. A fractional coloring $\gamma$ is called geometric if
\[
    \overline\gamma(S)=\overline\gamma(S')
\]
whenever $S,S'\subseteq V(G)$ are isometric as subsets of the plane. In other
words, the quantity $\overline\gamma(S)$ depends only on the congruence class of the
point configuration $S$, not on its particular location inside $G$.

The geometric fractional chromatic number of $G$, denoted $\chi_{gf}(G)$, is
the minimum weight of a geometric fractional coloring of $G$:
\[
    \chi_{gf}(G)
    =
    \min_{\gamma}
    \sum_{I\in\mathcal I(G)} \gamma(I) = \min_{\gamma} |\gamma|,
\]
where the minimum is taken over all geometric fractional colorings $\gamma$ of
$G$. Since geometric fractional colorings are fractional colorings with
additional linear constraints, we always have
\[
    \chi_f(G)\le \chi_{gf}(G).
\]

In fact, the original motivation behind the definition of the geometric fractional chromatic number is obtaining stronger bounds on $m_1(\mathbb{R}^2)$, since by \cite[Corollary 2]{ACMVZ} we also have

\[
    m_1(\mathbb{R}^2) \leq 1/\chi_{gf}(G)
\]
for any unit-distance graph $G$.

\section{The Augmented Graph and Proof of the Main Theorem}

\begin{figure}[H]
    \centering
    \includegraphics[width=0.5\linewidth]{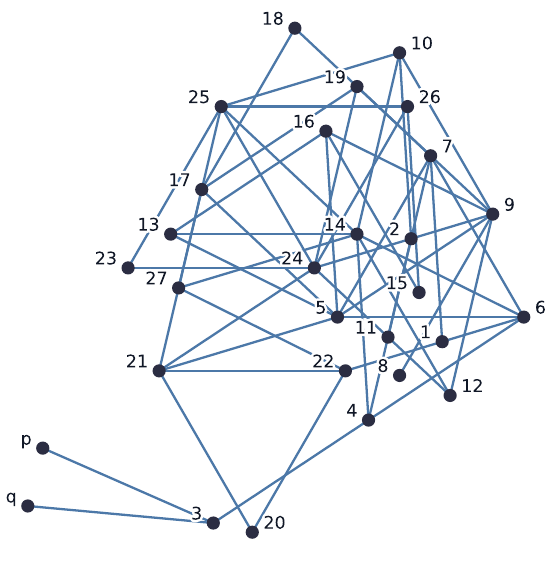}
    \caption{Our extension of the $G_{27}$ graph with vertices $p$ and $q$.}
    \label{fig:snail}
\end{figure}

Here we present the pair of vertices $p, q$, which achieve our $\chi_{gf}$ bound when added to $G_{27}$.
If we let $\omega_1 = 1/2 + i\sqrt{3}/2$ and $\omega_3 = 5/6 + i\sqrt{11}/6$, we can define the Moser lattice as\[L_{Moser} = \{a + b\omega_1 + c\omega_3 + d\omega_1\omega_3 : a,b,c,d\in\mathbb{Z}\}.\]
The entire $\G$ graph is a subset of $L_{Moser}$, and by \cite{DuczMoserLattice}, no graph in this lattice can achieve $\chi_{gf}(G) > 4$.

Vertex $p$ can be represented as

\[p=
3+\frac{17}{8}\omega_1-\frac78\omega_3+2\omega_1\omega_3
+\sqrt5\left(
-\frac14+\frac18\omega_1-\frac18\omega_3+\frac14\omega_1\omega_3
\right).\]

The point $q$ does not admit an equally simple description. It is given by

\[
q=
\frac{11}{4}
+\frac{13}{8}\omega_1
-\frac18\omega_3
+2\omega_1\omega_3
+
\frac{\eta}{8}
\left(-\omega_1+\omega_3+\omega_1\omega_3\right)
\]

where 

\[
\eta=
i\sqrt{\frac{415+79\sqrt{33}}8}
.
\]

The unit-distance graph $\Gp$ is defined by the points $V(\G) \cup \{p,q\}$, shown in Figure~\ref{fig:snail}. Its shape suggests the informal nickname ``Snail graph''.

\begin{figure}[H]
    \centering
    \includegraphics[width=0.5\linewidth]{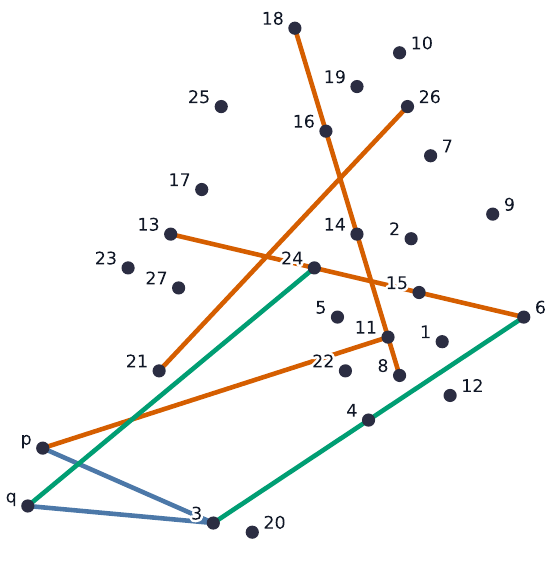}
    \caption{Congruences containing the vertices $p$ and $q$. The blue intervals both have length 1, the green intervals both have length 2, and the red intervals all have the same length.}
    \label{fig:snail_congs}
\end{figure}

The two new vertices only give rise to the trivial congruences $\{p\} \cong \{q\} \cong \{v_i\}$, and two nontrivial 2-element congruences, namely $\{p, v_{11}\} \cong \{v_{6}, v_{13}\} \cong \{v_{21}, v_{26}\} \cong \{v_8, v_{18}\}$ and $\{q, v_{24}\} \cong \{v_{6},v_{3}\}$. Both $p$ and $q$ are degree-1 vertices connected to $v_3$. The schematics of the new congruences is shown in Figure~\ref{fig:snail_congs}.

\begin{lemma}
 \[\chi_{gf}(\Gp) > 4.0007.\]
\end{lemma}

This is certified by a rational dual solution to the geometric fractional coloring number (GFCN) linear program of $\Gp$, which is given in the supplementary material \cite{supplementary}. The verification process is essentially the same as in \cite{MRVZ}, except that our rational certificate is not claimed to be optimal. Determining the exact value of $\chi_{gf}(\Gp)$ would require solving the GFCN dual linear program for $\Gp$ exactly over the rationals. This program has $16860$ variables and $498168$ constraints, and its exact solution was beyond our available computational resources.

\bigskip
Following the arguments in the proofs of Theorems 1 and 2 in \cite{MRVZ}, we can obtain finite ``blow-ups'' of our graph $\Gp$ with independence ratio arbitrarily close to $1/\chi_{gf}(\Gp)$, and in particular under $1/4$, proving our Theorem~\ref{thm1}:

\begin{theorem*}[Restatement of Theorem~\ref{thm1}]
There exists a finite unit-distance graph $G$ in the plane such that
$\alpha(G)/|V(G)|<1/4$.
\end{theorem*}

\section{Characterizing the optimal geometric fractional colorings of $G_{27}$}

The rest of this note explains how the vertices $p$ and $q$ were found.
We now characterize the optimal geometric fractional colorings of $G_{27}$, or more specifically we enumerate all extremal points of the polytope of optimal geometric fractional colorings. 

The first step is to write up the linear program defining
$\chi_{gf}(G_{27})$. A geometric fractional coloring is represented by a
nonnegative vector $x=(x_I)_{I\in\mathcal{I}(\G)}$, where $x_I$ is the weight
assigned to the independent set $I$. Following the notation of \cite{MRVZ}, the
normalization and congruence constraints are encoded by a vector $e$ and a matrix
$C$, so that the relevant primal linear program has the form
\[
    \min \langle \mathbf{1},x\rangle
\]
subject to
\[
    x\geq 0,\qquad \langle e,x\rangle=1,\qquad Cx=0.
\]
Here $\langle \mathbf{1},x\rangle=\sum_{I\in\mathcal{I}(\G)}x_I$ is the total
weight of the coloring.
The vector $e$ is a $0-1$ vector with dimension equal to $|\mathcal{I}(\G)|$, with ones corresponding to independent sets containing the first vertex of $V(\G)$. Each row of $C$ corresponds to a congruence constraint between two isometric subsets $S, S'$ of $V(\G)$, encoding $\overline\gamma(S) = \overline\gamma(S')$.

The certificate proving $\chi_{gf}(G_{27})=4$ is a dual vector $y$
such that the corresponding dual inequalities are all satisfied. For each
$I\in\mathcal{I}(\G)$, write $C_I$ for the column of $C$ indexed by $I$, and define
the dual slack
\[
    r_I = y^T C_I - 4 e_I + 1 .
\]
Thus $r_I\geq 0$ for every $I\in\mathcal{I}(\G)$. Let
\[
    \mathcal{I}_0=\{I\in\mathcal{I}(\G): r_I=0\}
\]
be the set of independent sets on which the dual inequality is tight.

Due to complementary slackness, if $x$ is an
optimal primal solution, then no weight can be placed on an independent set
whose dual inequality has positive slack:
\[
    r_I>0 \quad\Longrightarrow\quad x_I=0 .
\]
Consequently, every optimal geometric fractional coloring is supported on
$\mathcal{I}_0$. Surprisingly, $|\mathcal{I}_0| = 168$, in sharp contrast with $|\mathcal{I}(G_{27})|=182304$. That is, only a minuscule fraction of the independent sets of $G_{27}$ are allowed to participate in optimal geometric fractional colorings.

Conversely, every feasible point supported on $\mathcal{I}_0$
is automatically optimal. Indeed, suppose that
\[
    x\in \mathbb{R}_{\geq 0}^{\mathcal{I}_0},\qquad
    \langle e,x\rangle=1,\qquad
    Cx=0 .
\]
For each $I$ in the support of $x$, the equality $r_I=0$ gives
\[
    1 = 4e_I - y^T C_I .
\]
Multiplying by $x_I$ and summing over $I$ gives
\[
    \langle \mathbf{1},x\rangle
    =
    4\langle e,x\rangle - y^T Cx
    =
    4 .
\]

Let $C_0$ and $e_0$ be the
restrictions of $C$ and $e$ to the columns indexed by $\mathcal{I}_0$.
Thus we characterize the optimal face if we determine the restricted polytope:
\[
    F =
    \left\{
        x\in \mathbb{R}_{\geq 0}^{\mathcal{I}_0} :
        \langle e_0,x\rangle=1,\;
        C_0 x=0
    \right\}.
\]

The affine hull of $F$ has dimension $11$. It remains to enumerate the extremal points of $F$. For rational coefficients, we can get an exact solution for this task using the Double Description Method \cite{FukudaProdon, cddlib}.

It turns out that the polytope $F$ of optimal geometric fractional colorings of $\G$ has exactly $23$ extremal points; we denote this set by $E$. Two elements of $E$ are the integer colorings characterized in \cite{DuczMoserLattice}, and hence are $1$-fold colorings. For each of the fold numbers $3$, $6$, and $12$, there are two extremal colorings, while the largest fold number, $22$, is attained by three extremal colorings. Four representative examples are shown in Table~\ref{table:extremals}; the full list of $23$ extremal colorings is included in the supplementary material \cite{supplementary}.

\begin{table}[H]
\centering
\renewcommand{\arraystretch}{1.2}
\begin{tabular}{|c|c|}
\hline
\includegraphics[width=0.3\textwidth]{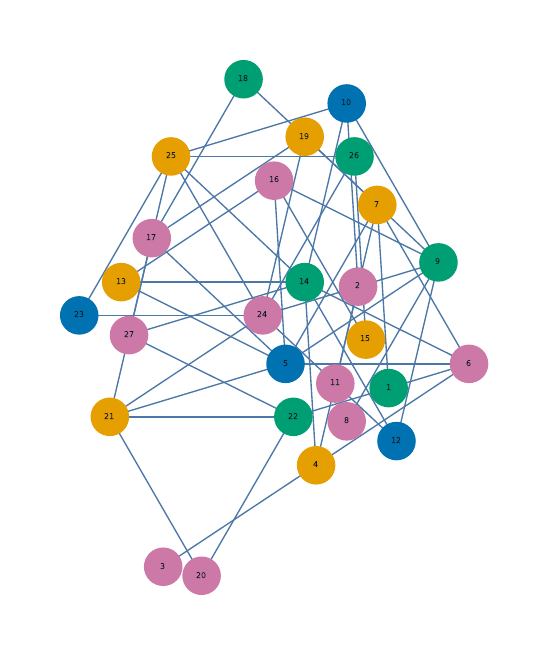} &
\includegraphics[width=0.3\textwidth]{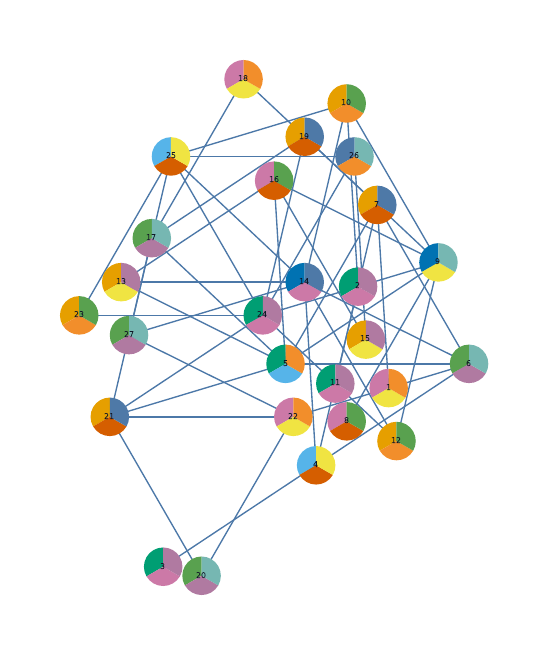} \\

\hline
\includegraphics[width=0.3\textwidth]{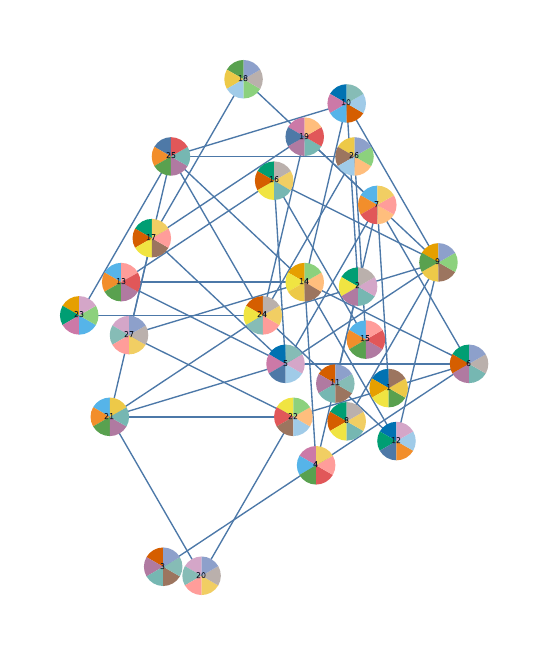} &
\includegraphics[width=0.3\textwidth]{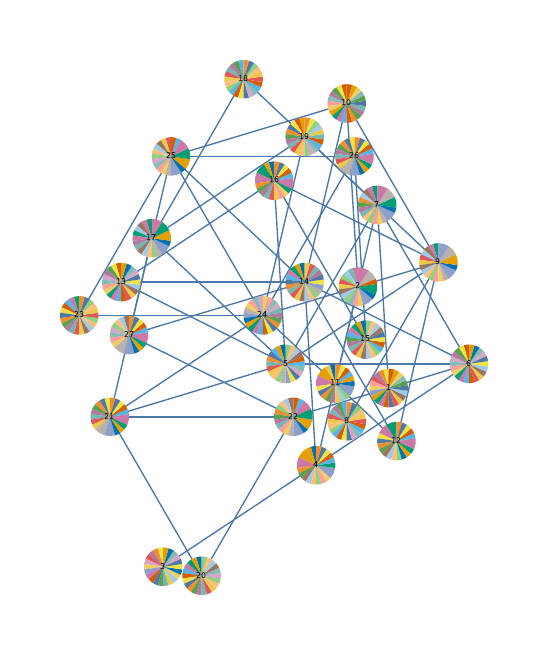} \\
\hline
\end{tabular}
\vspace{5pt}
\caption{Four examples from the set $E$ of $23$ extremal optimal geometric fractional colorings of $\G$. They are $1$-, $3$-, $6$-, and $22$-fold fractional colorings, respectively.}
\label{table:extremals}
\end{table}

\section{Finding suitable extra vertices}

With our characterization of all geometric fractional colorings of $\G$, we have a powerful tool for testing new vertex candidates.

\subsection{Generating candidate vertices}
Let $D$ be the set of all pairwise distances between the vertices of $G_{27}$. Our candidate vertex pool $C$ was generated by the following algorithm:

\begin{enumerate}
    \item Take any two vertices $v_i, v_j \in V(G_{27})$.
    \item Select two distances $d_1, d_2 \in D$.
    \item Let $c_1, c_2$ be the intersections of circles centered at $v_i, v_j$ with radii $d_1, d_2$ (if any exist). Add these to $C$.
\end{enumerate}

We then repeat these steps until a sufficiently large pool $C$ is found.

\subsection{Filtering good candidate pairs}

Due to the size of $V(G_{27})$, it is clearly infeasible to check all candidates in $C$ by evaluating their corresponding GFCN linear programs, since each candidate takes multiple minutes to evaluate. We therefore employ a prefiltering step on $C$, which only solves relatively small constraint satisfaction problems over integer variables.

The primal optima of the GFCN LP for the $G_{27}$ graph are all convex combinations of 23 extremal colorings. Since these colorings are rational with small denominators, we can interpret them as fractional colorings assigning $a$ colors out of a total of $b$ to each vertex of $G_{27}$, with $b = 4a$. If we let $c(x)$ denote the set of colors assigned to the vertex $x$, the GFCN constraints turn into:

\begingroup
\setlength{\abovedisplayskip}{3pt}
\setlength{\belowdisplayskip}{3pt}
\setlength{\abovedisplayshortskip}{3pt}
\setlength{\belowdisplayshortskip}{3pt}
\begin{align}
    \|x-y\|=1
    &\implies c(x)\cap c(y)=\emptyset, \label{cond:edge} \\
    S\cong S'
    &\implies
    \left|\bigcap_{x\in S} c(x)\right|
    =
    \left|\bigcap_{x\in S'} c(x)\right|. \label{cond:geometric}
\end{align}
\endgroup

For such a given rational coloring and a new candidate vertex $y$, we create a constraint satisfaction problem, where we must assign a set of exactly $a$ colors to $c(y)$, extending the previous coloring such that conditions $(1)$ and $(2)$ are met.
If this CSP is infeasible, then our particular coloring has no extension to the graph $G \cup \{y\}$. In our implementation, condition (2) was restricted to only $2$ and $3$-point congruences for efficiency.

Given candidate vertices $C$, and a set of rational geometric colorings $R$ of $G_{27}$, we can then efficiently check which vertices in $C$ obstruct which colorings of $R$. Hence, a promising approach to find a vertex $v \in C$ such that $\chi_{gf}(G \cup \{v\}) > 4$ is to choose $R$ to be equal to $E$, the 23-element set of extremal colorings. However, a vertex obstructing all extremal points $E$ of the optimal face $F$ does not imply that it obstructs all elements of $F$, and indeed the choice of $R=E$ only returned candidates with $\chi_{gf}(G \cup \{v\}) = 4$ during the allocated time budget.


Our chosen set of colorings $R$ therefore consisted of all of $E$, and additionally all possible $1:1$-ratio combinations of pairs of colorings from $E$, totalling 276 colorings checked for each vertex. A nearly 2-day search resulted in over $800$ vertex pairs $y_1, y_2 \in C$ that together block all $276$ colorings (we found no single vertex obstructing all of $R$). A large number of these vertex pairs appear to give the same or similar bounds numerically, with the best pair reaching a GFCN value of over $4.0009$, although the only pair rigorously verified is the one we present in this article.

If we have a blocking pair of vertices, then a rough estimate yields that one of its elements must obstruct at least 138 of the 276 colorings (disregarding the interaction between the newly added vertices). This property is exceedingly rare among vertices generated by the algorithm described in section 4.1. By randomly sampling elements of $C$, we estimate that less than $0.1\%$ of vertices have this property. Even obstructing just a single coloring is rare, with a probability close to $5\%$ per sampled vertex. For comparison, vertex $p$ blocks $273$ of the $276$ colorings, while $q$ blocks $252$.

\section{AI use declaration}

The authors used OpenAI's ChatGPT (GPT-5.5) during the preparation of this manuscript for language editing, grammar improvement, and formatting suggestions. The authors used OpenAI's ChatGPT and Codex (GPT-5.5) to assist with software development, including code drafting, debugging, and implementation of computational experiments. Codex was asked to help search for vertices of the optimal primal face of a linear program; the resulting code unexpectedly found all such vertices, not merely a sample. This observation influenced the subsequent research direction, in particular the decision to use the full characterization of the optimal face as a tool. All mathematical results and arguments were conceived, verified, and validated by the authors, who take full responsibility for the correctness of the work. Moreover, the supplementary verification script for this paper was written and verified entirely by hand.

\section{Acknowledgements}

 Á.~D. and D.~V. were supported by grant NKFIH-153165, and the Ministry of Innovation and Technology NRDI Office within the framework of the Artificial Intelligence National Laboratory (RRF-2.3.1-21-2022-00004). We thank Domonkos Czifra, Máté Matolcsi, Dömötör Pálvölgyi, and Pál Zsámboki for helpful discussions.

\end{document}